# Stochastic Ordering of Symmetric Core Event Probabilities With Respect To Row Margins For Multiple Hypergeometric Random Variables

**Bruce Levin**  BL6@COLUMBIA.EDU
*Department of Biostatistics*
*Mailman School of Public Health*
*Columbia University*
*New York, NY 10032, USA*

## Abstract

Suppose *X* is a frequency vector that follows a central multiple hypergeometric distribution, such as arises in random sampling of an *m*-category attribute from a finite population without replacement. We show that the probability that *X* satisfies some symmetrical—but otherwise arbitrary—interval constraints in each component decreases as the sample size increases or decreases away from one-half of the population size. The result is used to analyze the variance reduction in multinomial frequencies subject to arbitrary interval censoring.

**Keywords:** Interval censoring; multinomial distribution; multiple hypergeometric distribution; stochastic ordering; variance-reduction due to interval censoring.

## 1. Introduction and preliminaries

For any non-negative integer $n$, let $[n]$ denote the discrete interval of non-negative integers up to $n$, $[n] = \{0,1,...,n\}$. For given integers $m \geq 2$ and $n \geq 2$, let $\Omega_{m,n} = \{x \in [n]^m : x_1 + \cdots + x_m = n\}$ be the discrete simplex in *m* variables of order *n*. For given non-negative integers $s_1,...,s_m$ with total $t = s_1 + \cdots + s_m$, let *X* have the multiple hypergeometric distribution with row sum *n* and column sums $s = (s_1,...,s_m)$, say $X \sim H_m(n;s)$. The sample space for $H_m(n;s)$ is the intersection of the Cartesian product of discrete intervals $[s_j]$ with the discrete simplex: $\mathcal{H}_m(n;s) := \prod_{j=1}^{m}[s_j] \cap \Omega_{m,n}$.[1] As is well known, when sampling without replacement from an urn containing *t* balls of *m* different colors with $s_j$ balls of color $j$, $H_m(n;s)$ gives the distribution of the vector $X$ counting the number of balls of each color after $n$ draws.[2]

---

[1] The probability function is $P[X=x] = \prod_{j=1}^{m}\binom{s_j}{x_j} / \binom{t}{n} = \binom{n}{x}\binom{t-n}{s-x} / \binom{t}{n}$ where the first expression uses binomial coefficients and the second uses multinomial coefficients in the numerator, e.g., $\binom{n}{x} = n!/(x_1!\cdots x_n!)$. The terms "row sum" and "column sums" refer to the margins of a $2 \times m$ table containing the *X* frequencies in the first row and $s - X$ in the second row.

For $j=1,...,m$ and given integers $l_j, u_j \in [s_j]$ with $l_j \le u_j$ and $l_j + u_j = s_j$, consider the discrete intervals $B_j = \{l_j,...,u_j\}$ which are symmetrically placed around their symmetric centers $s_j/2$. We call the subset $B = B_{m,n}(l,u) = (B_1 \times \cdots \times B_m) \cap \Omega_{m,n}$ a *symmetric core* of $\mathcal{H}_m(n;s)$. We are interested in the probability of the event $[X \in B]$ and for different values of $n$, we'll write $P_n[X \in B]$ when $X \sim H_m(n;s)$. In this note we prove the following stochastic ordering property for symmetric core event probabilities with respect to the row sums.

<u>Theorem 1</u>. *For any symmetric core subset $B$ of $\mathcal{H}_m(n;s)$ and any two sample sizes $n$ and $n'$ with $\lceil t/2 \rceil \le n < n' \le t$,*

$$P_n[X \in B] \ge P_{n'}[X \in B]. \qquad (1)$$

*The inequality also holds for $0 \le n' < n \le \lfloor t/2 \rfloor$.*

Theorem 1 seems intuitively obvious when one considers special cases such as the "bet on books" discussed in Levin (1983). A shuffled deck of cards is evenly divided between two players. What would be fair odds to bet that neither player has any "books" of any rank?[3] To answer this requires the probability $P_{26}[X \in B] \approx 0.231453$ where $X \sim H_{13}(26;4,...,4)$ and $B_j = \{1,2,3\}$ for each $j$. (See Levin, 1983 for the calculation.) Intuition suggests that this probability would only decrease if the deck were unevenly divided.[4] The essence of Theorem 1 is that the stochastic ordering displayed in (1) holds even if books are permitted for some ranks, $B_j = \{0,...,4\}$, or the cards in some ranks are required to be evenly shared, $B_j = \{2\}$.

The restriction to symmetric core subsets is essential as (1) does not generally hold for asymmetrical discrete interval constraints, $R_j$ say, on the cell frequencies. For example, if $X$ has the simple hypergeometric distribution $X \sim H_2(n;4,6)$, with asymmetrical interval constraints $R_1 = \{0,...,3\}$ and $R_2 = \{3,...,6\}$, then $P_5[X \in (R_1 \times R_2) \cap \Omega_{2,5}] \approx 0.74$ but $P_6[X \in (R_1 \times R_2) \cap \Omega_{2,6}] \approx 0.93$, violating (1).

---

[2] As is also well known, $H_m(n;s)$ is the conditional distribution of $m$ independent binomial random variables $Y_j \sim B(s_j, p)$ given fixed sum $Y_1 + \cdots + Y_m = n$ for any $0 < p < 1$ and it is also the conditional distribution of two independent multinomial random variables, $Y_1 \sim M_m(n, P)$ and $Y_2 \sim M_m(t-n, P)$ given fixed sum $Y_1 + Y_2 = s$ for any $P$ in the continuous simplex $\Delta_m = \{p \in \mathfrak{R}^m_{\ge 0} : p_1 + \cdots + p_m = 1\}$.

[3] A *book* comprises all four cards of a given rank.

[4] Even equalling 0 if $n > 39$. When $n = 27$ or 25, $P_n[X \in B] \approx 0.225406$.

Theorem 1 appears not to have been discussed in the literature. Block et al. (1985) discuss a property they call *negative dependence through stochastic ordering* (NDS) and show that the multiple hypergeometric distribution is NDS, which implies that $E[g(Y)|\sum_j Y_j = n]$ is non-increasing in *n* for any non-decreasing Borel measurable function *g* for which the conditional expectations exist. With independent $Y_j \sim Bin(s_j, \frac{1}{2})$, it would be natural to consider taking $g(Y) = I[Y \in B]$, the conditional expectation of which is $P[X \in B]$, but *g* both increases and decreases in any coordinate for which $l_j > 0$ and $u_j < s_j$, so an immediate application of NDS to prove Theorem 1 is not apparent. Klenke and Mattner (2010) discuss stochastic ordering for univariate hypergeometric tail probabilities with respect to one of the margins, but their results do not generalize to probabilities of symmetric core subsets which are *differences* of tail probabilities.

After proving Theorem 1 in Section 2, we mention an application in the theory of truncated or interval-censored multinomial frequencies in Section 3. The connections are detailed in Levin (2023).

## 2. Proof of Theorem 1

With $n \geq \lceil t/2 \rceil$, it clearly suffices to consider $n' = n+1$ and extend the result to larger $n'$ by induction. Without loss of generality we may also assume that it is numerically possible to realize either endpoint cell frequency $l_j$ or $u_j$ for each *j*.[5] Given a sample $X = x$ from $H_m(n;s)$, consider drawing one additional ball from the urn, in which after *n* draws there remains $t-n$ balls with $s_j - x_j$ of color *j*. Let *C* denote the color of the additional ball and let $Y = (I[C=1],...,I[C=m])$ be the frequency count vector of the additional draw with $Y \sim H_m(1; s-x)$.[6] We'll write $X' = X + Y$ for the final frequencies.

Among the several ways the event $[X \in B]$ can occur if one "sticks" with the sample of size *n*, only certain types of event affect the difference $P_n[X \in B] - P_{n+1}[X' \in B]$. One type is of the form $E_J = [X \in B, X_j = u_j \text{ for } j \in J \text{ and } X_j < u_j \text{ for } j \notin J]$ for subsets $J \subseteq \{1,...,m\}$ of size at least 1.[7] In words, in event $E_J$, one or more cell frequencies *j* have reached their maximum quota $u_j$, namely, precisely those cells in *J*, with all other cells satisfying their minimum and maximum quotas, the latter strictly. For such events, after an additional ball is drawn, a "win", i.e., $[X' \in B]$, requires *C* to avoid all the cells in *J*. For convenience, we'll write $E_J = [X_J = u_J, X_{-J} < u_{-J}, X_{-J} \in B_{-J}]$, where $X_J$ is the sub-vector of *X* including only coordinates $j \in J$ (and similarly for $u_J$); $X_{-J}$ is the sub-vector of *X*

---
[5] Else we could move the endpoints closer to $s_j/2$ without changing the value of $P_n[X \in B]$.

[6] Since $n=1$, *Y* is also multinomial with $Y \sim M_m(1, (s-x)/(t-n))$.

[7] The event $E_J$ for $J = \{1,...,m\}$ is possible if $u_1 + \cdots + u_m = n$ but impossible otherwise.

−3−

including only coordinates $j \notin J$ (and similarly for $u_{-J}$); and where $B_{-J}$ is the intersection of the sub-Cartesian product of discrete intervals including $B_j$ only for $j \notin J$ with the discrete simplex $\Omega_{m^-, n^-}$ for $m^- = m - |J|$ and $n^- = n - \sum_{j \in J} u_j$.

Another type of event, disjoint from the first, is of the form $F_j = [X_j = l_j - 1, X_k \in B_k \text{ for } k \neq j]$ $= [X_j = l_j - 1, X_{-j} \in B_{-j}]$ for singleton sets $\{j\}$, i.e., where one and only one cell frequency falls short of its minimum quota by one ball, with all other cells satisfying their quotas. Clearly, no more than a single cell can be one ball short of its minimum quota in order for a single additional draw to produce a win.[8] For $F_j$ events, a win with an additional draw requires $C=j$, but without the additional draw, $[X \notin B]$. Only the two types of events $E_J$ and $F_j$ affect the difference $P_n[X \in B] - P_{n+1}[X' \in B]$ because any other $X$ leading to $[X \in B]$ after $n$ draws would also lead to $[X' \in B]$, irrespective of $Y$. Therefore, the point probabilities for such other winning $X$, say $X = x^*$, cancel in the difference because the individual terms match, i.e., $P_{n+1}[X = x^*, X' \in B] = \sum_{\text{any } y} P_n[X = x^*, C = y] = P_n[X = x^*]$.[9]

We may now write,

$$P_n[X \in B] - P_{n+1}[X' \in B] = \sum_J P_n[E_J] - \{\sum_J P[E_J, C \notin J] + \sum_j P_n[F_j, C = j]\}$$
$$= \sum_J P_n[E_J, C \in J] - \sum_j P_n[F_j, C = j], \tag{2}$$

where the first sum in the final expression is over all non-empty subsets $J \subseteq \{1,...,m\}$ and the second sum is over all singleton indices $j \in \{1,...,m\}$. For $j = 1,...,m$, let $\Im_j$ denote the subset of all $x \in B$ such that $x_j = u_j$. It is easy to see that $\Im_j$ is the disjoint union of $E_j$ with $\bigcup_{j' \neq j} E_{\{j,j'\}}$ and $\bigcup_{j,j',j''} E_{\{j,j',j''\}}$ (for distinct $j, j', j''$) and so forth. Also note that $\Im_j = [X_j = u_j, X_{-j} \in B_{-j}]$, with no restrictions on the remaining sub-vector $X_{-j}$ besides satisfying the remaining quotas in $B_{-j}$ and summing to the reduced sample size $n - u_j$. *We now claim that the first sum in (2) can be re-expressed solely in terms of singleton indices as*

---

[8] For convenience we drop the braces around singleton subscripts, writing $E_j$ for $E_{\{j\}}$ and $F_j$ for $F_{\{j\}}$. In $F_j$, $X_{-j}$ has sample size $n - l_j + 1$, so $X_{-j} \in \Omega_{m^-, n^-}$ where now $m^- = m - 1$ and $n^- = n - l_j + 1$. The context will always make clear what values $m^-$ and $n^-$ assume in $B_{-J}$ for $E_J$ or $B_{-j}$ for $F_j$. We may also assume without loss of generality that $u_1 + \cdots + u_m > n$ for if not, we must have $X = u$, but then $F_j = \phi$ for each $j$.

[9] Obviously, events wherein some cells exceed their upper quotas enter neither term of $P_n[X \in B] - P_{n+1}[X' \in B]$.

$$\sum_J P_n[E_J, C \in J] = \sum_j P_n[X_j = u_j, X_{-j} \in B_{-j}, C = j]. \qquad (3)$$

This is because for each subset $J = \{j, j', j'',...\}$ with two or more elements, we write $P_n[E_J, C \in J] = P_n[E_J, C = j] + P_n[E_J, C = j'] + P_n[E_J, C = j''] + \cdots$ and then collect terms with common $j$, as follows.

$$\sum_J P_n[E_J, C \in J] = \sum_j P_n[E_j, C = j] + \left\{ \sum_{j' \neq j} P_n[E_{\{j,j'\}}, C = j] + \sum_{j' \neq j} P_n[E_{\{j,j'\}}, C = j'] \right\}$$

$$+ \left\{ \sum_{\substack{distinct \\ j,j',j''}} P_n[E_{\{j,j',j''\}}, C = j] + \sum_{\substack{distinct \\ j,j',j''}} P_n[E_{\{j,j',j''\}}, C = j'] + \sum_{\substack{distinct \\ j,j',j''}} P_n[E_{\{j,j',j''\}}, C = j''] \right\} + \cdots$$

$$= \sum_j \left\{ P_n[E_j, C = j] + \sum_{j' \neq j} P_n[E_{\{j,j'\}}, C = j] + \sum_{\substack{distinct \\ j,j',j''}} P_n[E_{\{j,j',j''\}}, C = j] + \cdots \right\}$$

$$= \sum_j P_n[E_j \cup \bigcup_{j' \neq j} E_{\{j,j'\}} \cup \bigcup_{\substack{distinct \\ j,j',j''}} E_{\{j,j',j''\}} \cup \cdots, C = j]$$

$$= \sum_j P_n[\Im_j, C = j] = \sum_j P_n[X_j = u_j, X_{-j} \in B_{-j}, C = j].$$

Thus (2) becomes

$$P_n[X \in B] - P_{n+1}[X' \in B] = \sum_j \left\{ P_n[X_j = u_j, X_{-j} \in B_{-j}, C = j] - P_n[X_j = l_j - 1, X_{-j} \in B_{-j}, C = j] \right\} \qquad (4)$$

since by definition of $F_j$, $P_n[F_j, C = j] = P_n[X_j = l_j - 1, X_{-j} \in B_{-j}, C = j]$. *It will therefore suffice to show that the individual terms inside the braces in (4) are non-negative for each j.*

Going forward, we may assume without loss of generality that at least one lower limit $l_j$ is positive, for if not, then $F_j = \phi$ for each $j$, so $\sum_j P[F_j, C = j] = 0$ and $\delta_n \geq 0$. We may also assume that $B \neq [s_1] \times \cdots \times [s_m] \cap \Omega_{m,n}$, for obviously both outcome probabilities would equal 1 and $\delta_n$ would equal 0. Now for any given $j$, write

$$P_n[X_j = u_j, X_{-j} \in B_{-j}, C = j] = P_n[X_j = u_j] \cdot P_n[X_{-j} \in B_{-j} \mid X_j = u_j] \cdot P_n[C = j \mid X_j = u_j, X_{-j} \in B_{-j}]$$

and

$$P_n[X_j = l_j - 1, X_{-j} \in B_{-j}, C = j] = P_n[X_j = l_j - 1] \cdot P_n[X_{-j} \in B_{-j} \mid X_j = l_j - 1] \cdot P_n[C = j \mid X_j = l_j - 1, X_{-j} \in B_{-j}].$$

It follows that to show each term in (4) is non-negative is equivalent to showing

$$\frac{P_n[X_{-j} \in B_{-j} \mid X_j = u_j]}{P_n[X_{-j} \in B_{-j} \mid X_j = l_j - 1]} \geq \frac{P_n[X_j = l_j - 1]}{P_n[X_j = u_j]} \cdot \frac{P[C = j \mid X_j = l_j - 1, X_{-j} \in B_{-j}]}{P[C = j \mid X_j = u_j, X_{-j} \in B_{-j}]}. \tag{5}$$

We proceed to prove (5).

Note that given $X = x$, $C$ is independent of the individual elements of $X_{-j}$ because the conditional probability $P[C = j \mid X = x] = (s_j - x_j)/(t - n)$ depends on $x$ only through $x_j$. Furthermore, the marginal distribution of $X_j$ is hypergeometric, $X_j \sim H_2(n; s_j, t - s_j)$. Therefore, the right-hand side of (5) equals

$$\frac{\binom{s_j}{l_j - 1}\binom{t - s_j}{n - l_j + 1}}{\binom{s_j}{u_j}\binom{t - s_j}{n - u_j}} \cdot \frac{s_j - l_j + 1}{s_j - u_j}.$$

However, temporarily suppressing the subscript $j$ in the above expression, we find

$$\frac{\binom{s}{l-1}}{\binom{s}{u}} \cdot \frac{s - l + 1}{s - u} = \frac{u!(s - u - 1)!}{(l - 1)!(s - l)!} = \frac{\binom{s-1}{s-l}}{\binom{s-1}{u}} = 1 \quad \text{since } l + u = s.$$

Thus, the right-hand side of (5) equals $\binom{t-s}{n-l+1} / \binom{t-s}{n-u}$ and we claim this is $<1$ for $n \geq \lceil t/2 \rceil$.

<u>Proof of claim</u>. Note that because $u > l - 1$, we have (i) $n - l + 1 > n - u$, and by assumption, $n \geq \lceil t/2 \rceil$ so that $2n + 1 > t$; thus we have (ii) $n - l + 1 > t - n - l = t - s - n + u$. But either $n - u$ or $t - s - n + u$ is at least $(t - s)/2$, whence $n - l + 1 > max \geq (t - s)/2 > (t - s - 1)/2$ where $max = (n - u) \vee (t - s - n + u)$. Therefore, because $\binom{t-s}{x+1} < \binom{t-s}{x}$ if and only if $x > (t - s - 1)/2$, it follows that $\binom{t-s}{x'} < \binom{t-s}{x}$ for $x' > x > (t - s - 1)/2$, whence $\binom{t-s}{n-l+1} < \binom{t-s}{max} = \binom{t-s}{n-u} = \binom{t-s}{t-s-n+u}$. Thus, $\binom{t-s}{n-l+1} / \binom{t-s}{n-u} < 1$ as claimed. □

On the left-hand side of (5), we note that given $X_j = u_j$, $X_{-j}$ has the same conditional distribution as $W_j \sim H_{m-1}(n - u_j; s_{-j})$ and given $X_j = l_j - 1$, $X_{-j}$ has the same conditional distribution as $W'_j \sim H_{m-1}(n - l_j + 1; s_{-j})$. Notice that by symmetry of the core subset $B_{-j}$, $P_{n-u_j}[W_j \in B_{-j}] = P_{t-s_j-n+u_j}[W_j \in B_{-j}]$, these being the same probabilities from the urn experiment (among balls withdrawn or left in the urn). But since $n - l_j + 1 > max \geq (t - s_j)/2$ as above, we can argue by induction on $m$ that

–6–

$$P_n[X_{-j} \in B_{-j} | X_j = u_j] = P_{max}[W_j \in B_{-j}] \geq P_{n-l_j+1}[W'_j \in B_{-j}] = P_n[X_{-j} \in B_{-j} | X_j = l_j - 1]$$

so that the left-hand side of (5) is not less than 1 whereas the right-hand side of (5) is less than 1 as previously shown. Thus (5) holds.

It remains only to start the induction at $m=2$. Suppose $j=1$ (the argument for $j=2$ is analogous). By the same logic leading to (5), it suffices just to show $P_n[X_2 \in B_2 | X_1 = u_1] \geq P_n[X_2 \in B_2 | X_1 = l_1 - 1]$. Now, if $X_1 = l_1 - 1$, in order for $X_2 \in B_2$ to occur, we require $u_2 \geq X_2 = n - X_1 = n - l_1 + 1$, i.e., we require $n < l_1 + u_2$ for $X_2 \in B_2$ to be possible. If this is not the case, then obviously $P_n[X_2 \in B_2 | X_1 = u_1] \geq 0 = P_n[X_2 \in B_2 | X_1 = l_1 - 1]$. So assume that

$$n < l_1 + u_2. \tag{6}$$

We now show that if $X_1 = u_1$ then $X_2 \in B_2$. First, from (6), $n < l_1 + u_2 \leq u_1 + u_2$ so that $X_2 = n - X_1 = n - u_1 \leq u_2$. Second, $X_2 \geq l_2$, for if not, then $X_2 = n - u_1 < l_2$, which implies that $n < l_2 + u_1 = (s_2 - u_2) + (s_1 - l_1) < s_1 + s_2 - n$ by (6). But then $2n < t$, which contradicts the assumption $n \geq \lceil t/2 \rceil \geq t/2$. Thus, whenever $B$ and $n$ are such that $P_n[X_2 \in B_2 | X_1 = l_1 - 1] > 0$, we also have $P_n[X_2 \in B_2 | X_1 = u_1] = 1 \geq P_n[X_2 \in B_2 | X_1 = l_1 - 1]$. This completes the initial step of the induction.

The proof of (1) for the case $0 \leq n' < n \leq \lfloor t/2 \rfloor$ follows by symmetry upon applying (1) to the complementary variable $s - x$ with sample size $t - n$. This concludes the proof of Theorem 1. □

The following lemma gives a condition which is equivalent to the stochastic ordering of Theorem 1. The lemma follows from an application of Bayes' rule and is in fact the basis for the efficient computation of multiple hypergeometric probabilities with arbitrary discrete interval constraints on the cell frequencies (i.e., not necessarily symmetric core events) as well as for other discrete distributions. See Levin (1981, 1983, 1992, and 2014) for details. The use of $Y$ and $W$ below differs from that used in the proof of Theorem 1.

<u>Lemma 1</u>. *Let $Y_1, ..., Y_m$ be independent binomial random variables with $Y_j \sim Bin(s_j, \frac{1}{2})$ and let $W_j$ be $Y_j$ restricted to the interval $R_j = \{l_j, ..., u_j\}$ for any given $l = (l_1, ..., l_m)$ and $u = (u_1, ..., u_m)$ with $0 \leq l_j \leq u_j \leq s_j$. Let $Y = Y_1 + \cdots + Y_m \sim Bin(t, \frac{1}{2})$ and let $W = W_1 + \cdots + W_m$, whose distribution is the convolution of the m statistically independent truncated binomial random variables. Then for the multiple hypergeometric random variable $X \sim H_m(n;s)$ and event $R(l,u) = (R_1 \times \cdots \times R_m) \cap \Omega_{m,n}$,*

$$P[X \in R(l,u)] = P[W = n] \prod_{j=1}^{m} P[Y_j \in R_j] / P[Y = n]. \qquad \square$$

Applying this result to a symmetric core subset $B$ of $\mathcal{H}_m(n;s)$, we see that the ordering of Theorem 1 holds if and only if the ratio of successive point probabilities for the convolution $W$ is not less than that of the untruncated binomial sum $Y$. This yields the following result.

<u>Corollary to Theorem 1</u>. *The convolution of truncated binomial random variables W satisfies*

$$\frac{P[W=n]}{P[W=n+1]} \geq \frac{P[Y=n]}{P[Y=n+1]} = \frac{\binom{t}{n}}{\binom{t}{n+1}} = \frac{n+1}{t-n} \text{ for } n \text{ satisfying } \lceil t/2 \rceil \leq n < t. \qquad \square$$

As is well known, binomial random variables are discrete log-concave as are truncated binomials, and so is the convolution of the truncated binomials. See, e.g., Saumard and Wellner (2014). It can also be shown that the variance of $W$ is less than the variance of $Y$ for non-trivial $B$ (in fact, for arbitrary non-trivial rectangular events $R$). For symmetric core subsets, we can think of the prototypical continuous model of two normal pdf's $f_W(x)$ and $f_Y(x)$ with the same mean $\mu$ but different variances, $\sigma_Y^2 > \sigma_W^2$. For that model, it is easy to show that $f_W(x)/f_W(x') > f_Y(x)/f_Y(x')$ for any $x' > x \geq \mu$. The corollary proves that this heuristic model holds as well for the discrete prototype case of symmetrical binomial random variables with symmetrical truncation about their means.

## 3. Discussion

A need for Theorem 1 arose in an investigation of the variance reduction that occurs in linear combinations of multinomial cell frequencies under discrete interval constraints. Given a probability vector $p \in \Delta_{m,n}$, suppose a multinomial random variable $X \sim Mult_m(n,p)$ is observed under interval constraints $X \in R$ where $R = R(l,u) = (R_1 \times \cdots \times R_m) \cap \Omega_{m,n}$ for *arbitrary* discrete intervals $R_j = \{l_j,...,u_j\}$ with $0 \leq l_j \leq u_j \leq n$. Under the constraints, the expected value of $X$ will generally not equal the unconstrained mean $np$ but will have a different conditional mean value, say $\mu = \mu(p,R) = E_p[X \mid X \in R]$. Furthermore, the variance of a given linear combination $c^T X$ will typically be strictly less than $c^T \{Diag[\mu] - \mu\mu^T/n\}c$, which is the variance of $c^T Y$ where $Y \sim Mult_m(n, \mu/n)$.[10] However, depending on $R$ and $c$, there may be no variance reduction at all. For example, if $m=3$, $l_1 = l_2 = l_3 = 0$, and $u_1 = u_2 = n$, the linear combination $c^T X$ with $c = (-p_2, p_1, 0)^T$ has exactly the same variance as $c^T Y$ for any $u_3 < n$. Under these circumstances it is reasonable to ask whether there is always some positive variance reduction in individual components $X_j$ of $X$, even when there is no interval constraint on some $X_j$, such as in the example for $j=1,2$. The answer to this question is yes, so long as $R$ imposes non-trivial constraints on at least one component. This is perhaps unsurprising given the negative dependence between components of *X*, but how to sculpt this fact into a rigorous argument is far from obvious for $m > 2$ and one would be well-justified in demanding a valid proof of the assertion. Levin (2023) gives such a proof wherein Theorem 1 plays the clinching role in establishing the result. Levin (2023) also fully explains when linear combinations of constrained frequencies exhibit no variance reduction.

---

[10] Note that here we're comparing variances under the *truncated* multinomial distribution versus an *untruncated* multinomial *with the same expectation* $\mu$ as that of the truncated distribution.